\documentclass[12pt,english]{article}
\usepackage[margin=0.75in]{geometry}
\usepackage[utf8]{inputenc}
\usepackage[english]{babel}
\usepackage{amsthm}
\usepackage{graphicx}
\usepackage{url}
\usepackage{hyperref}
\usepackage{amsmath}
\usepackage{listings}
\usepackage{amsfonts}
\usepackage{xcolor}
\usepackage{float}
\newtheorem{theorem}{Theorem}

\def\XXint#1#2#3{{\setbox0=\hbox{$#1{#2#3}{\int}$}
     \vcenter{\hbox{$#2#3$}}\kern-.5\wd0}}

\begin{document}

\title{On the series expansion of the secondary zeta function about $s=1$ and its coefficients}

\author{Artur Kawalec}

\date{}
\maketitle

\begin{abstract}
The secondary zeta function is defined as a generalized zeta series over the imaginary parts of non-trivial zeros assuming (RH). This function admits Laurent series expansion at the double pole at $s=1$. In this article, we derive a new formula for the expansion coefficients of the regular part, which is similar to the Stieltjes constants formula for the Riemann zeta function. We also numerically verify and compute the new formula to high precision for several test cases. Lastly, we also apply the Brent's (BPT) Theorem for improving convergence of the main formula.
\end{abstract}

\section{Introduction}
Let $\gamma_n$ denote the positive imaginary parts of non-trivial zeros of the Riemann zeta function, written as $\rho=\tfrac{1}{2}+i\gamma_n$ on the critical line assuming (RH), where the first few values are $\gamma_1=14.13472514\ldots$, $\gamma_2=21.02203963\ldots$, $\gamma_3=25.01085758\ldots$, and so on. The secondary zeta function is defined as a generalized zeta series over these imaginary parts of non-trivial zeros as

\begin{equation}\label{eq:1}
Z(s)=\sum_{\gamma>0}\frac{1}{\gamma^s}
\end{equation}
which converges absolutely for $\Re(s)>1$. The secondary zeta function extends to $\mathbb{C}\backslash\{1,-1,-3,-5,\ldots\}$ with a double pole at $s=1$ and simple poles at negative odd integers. The Laurent series expansion about $s=1$ is given by

\begin{equation}\label{eq:1}
Z(s) = \frac{1}{2\pi(s-1)^2} - \frac{\log(2\pi)}{2\pi(s-1)} + \sum_{n=0}^{\infty} C_n \frac{(s-1)^n}{n!}
\end{equation}
as derived by Ivi$\acute{c}$ [14][6]. This series has radius of convergence $R=2$ due to being limited by the next closest pole at $s=-1$ from the center at $s=1$. In this article, we derive a new formula for these expansion coefficients as

\begin{theorem}
\label{thm:fta}
\begin{equation}\label{eq:1}
C_n =\lim_{T\to\infty} (-1)^n \Bigg\{\sum_{0<\gamma<T}\frac{\log^n(\gamma)}{\gamma}-\frac{1}{2\pi(n+1)(n+2)}\log^{n+1}(T)\log\left(\frac{T^{n+1}}{(2\pi)^{n+2}}\right)\Bigg\}
\end{equation}
for $n\geq 0$.
\end{theorem}

We first review some relevant sources on this subject in the References Section. Chakravarty first wrote down the principal part of (2), and coined the term secondary zeta function [7][8][9]. And Ivi$\acute{c}$ derived the main Laurent series expansion with regular part in (2). Hassani in [12][13] proved the limit (3) for the $n=0$ case, and who also first computed the constant $C_0$ to $5$ decimal places. More developments by Brent (BPT) [3][4][5], who developed better error bounds for the sums over non-trivial zeros, and who improved computation of Hassani constant $C_0$ to $19$ decimal places [4, p.64]. A key development by  Arias De Reyna [2][15] is the (ADR) algorithm for high precision computation of $Z(s)$ for any complex $s$, which we used to compute $C_n$ to high precision in [16] and published a reference table of $C_n$ coefficients to $50$ accurate digits (see also Table $1$ below for $C_n$ coefficients). Also in [17], we published an unprocessed list of over $100$ coefficients, which took over $1$ month to compute on a fast workstation.  We believe that many of these coefficients should be accurate to more than $100$ digits after the decimal place, but more analysis is still needed.

In this article, based on these results, we derive a general limit formula (3) for $C_n$ for $n\geq 0$ and connect it to the series expansion (2).

\section{Proof of Theorem 1}
To show that, let us recap some basic definitions. The non-trivial zero counting function up to height $T$ (assuming $T\neq \gamma$) is defined by

\begin{equation}\label{eq:1}
N(T)=\sum_{0<\gamma<T} 1
\end{equation}
or by the average value

\begin{equation}\label{eq:1}\nonumber
N(T)=\frac{1}{2}\left[\sum_{0<\gamma<T}1+\sum_{0<\gamma\leq T}1\right]
\end{equation}
when taking into account $\gamma=T$, where there is a $\frac{1}{2}$ jump. But for the rest of this article, we consider the sums $0<\gamma<T$ for simplicity. The Riemann-von Mangoldt asymptotic formula is

\begin{equation}\label{eq:1}
N(T)=L(T)+Q(T)
\end{equation}
where the dominant term is

\begin{equation}\label{eq:1}
L(T)=\frac{7}{8}+\frac{T}{2\pi}\log\left(\frac{T}{2\pi}\right)-\frac{T}{2\pi}
\end{equation}
and the lower order term
\begin{equation}\label{eq:1}
Q(T)=S(T)+f(T)
\end{equation}
where

\begin{equation}\label{eq:1}
S(T)=\frac{1}{\pi}\arg\zeta(\tfrac{1}{2}+iT)
\end{equation}
defined as in [14]. One also has the bounds
\begin{equation}\label{eq:1}
S(T)=O(\log T), \quad f(T)=O(\frac{1}{T}), \quad Q(T)=O(\log T).
\end{equation}
Also, the $m^{th}$ derivatives of (1) are given by
\begin{equation}\label{eq:1}
(-1)^mZ^{(m)}(s)=\sum_{\gamma>0}\frac{\log^m(\gamma)}{\gamma^s}
\end{equation}
for $\Re(s)>1$. So we proceed by modifying the main derivation by Brent (BPT) [4, p.61] by inserting $\log^m(t)$ factor (setting $s=1$), hence we obtain

\begin{align}
\sum_{0<\gamma<T}\frac{\log^m(\gamma)}{\gamma} &= \int_{1}^{T} \frac{\log^m(t)}{t}dN(t)
= \int_{1}^{T} \frac{\log^m(t)}{t}dL(t) + \int_{1}^{T} \frac{\log^m(t)}{t}dQ(t) \notag \\[1.2em]
&= \frac{1}{2\pi} \int_{1}^{T} \frac{\log^m(t)\log(t/2\pi)}{t}\, dt
+ \left[ \frac{\log^m(t)}{t}Q(t) \right]_{1}^{T} - \int_{1}^{T} \frac{d}{dt}\left[\frac{\log^{m}(t)}{t}\right]Q(t)\, dt \notag \\[1.2em]
&= A(T) + \frac{\log^m(T)}{T}Q(T)+B_m+\int_{1}^{T}\frac{\log^m(t)-m\log^{m-1}(t)}{t^2}Q(t)\, dt
\end{align}
Here we note there is a subtle detail concerning the lower integration variable constant $B_m$ defined by

\begin{equation}\label{eq:1}
B_m=
\begin{cases}
-Q(1) & \text{if } m = 0 \\
0   & \text{if } m > 0
\end{cases}
\end{equation}
which only appears for $m=0$ case. The reason is that it is assumed the convention that $0^0=1$, which implies $[\log(1)]^0=1$, thus yielding a nonzero term, otherwise for $m>1$, we have $[\log(1)]^m=0^m=0$. Next, following this result, this integral identity can be generated as

\begin{equation}\label{eq:1}
\begin{aligned}
A(T)&=\frac{1}{2\pi} \int_{1}^{T} \frac{\log^m(t)\log(t/2\pi)}{t}\, dt\\[1.2em]
&=\frac{1}{2\pi(m+1)(m+2)}\log^{m+1}(T)\log\left(\frac{T^{m+1}}{(2\pi)^{m+2}}\right)
\end{aligned}
\end{equation}
by a repeated integration by parts for $m\geq 0$. Then, it is readily seen that the limit

\begin{equation}\label{eq:1}
\sum_{0<\gamma<T}\frac{\log^m(\gamma)}{\gamma}-A(T)=(-1)^mC_m+O\left(\frac{\log^{m+1}(T)}{T}\right)
\end{equation}
exists as $T\to\infty$ and converges to a constant

\begin{equation}\label{eq:1}
(-1)^mC_m=B_m+\int_{1}^{\infty}\frac{\log^m(t)-m\log^{m-1}(t)}{t^2}Q(t)\, dt
\end{equation}
for $m\geq 0$, since $Q(T)=\log(T)$. Now the next step is to connect this constant to Laurent Series expansion (2). Ivi$\acute{c}$ in [14] essentially derived this result, so we reproduce it as follows:

\begin{equation}\label{eq:1}
\begin{aligned}
Z(s)
&= \int_{1}^{\infty}\frac{1}{t^{s}}\, dN(t) = \int_{1}^{\infty} \frac{1}{t^{s}} \Biggl[
  \frac{1}{2\pi} \log\left(\frac{t}{2\pi}\right)\Biggr]dt\,  +  \int_{1}^{\infty} \frac{1}{t^{s}}dQ(t)
\\[1.2em]
&= \frac{1}{2\pi} \Biggl[
  t^{1-s} \log\left(\frac{t}{2\pi}\right) \Biggr]_1^\infty
  - \frac{1}{2\pi} \int_1^\infty \frac{t^{1-s}}{1-s} \cdot \frac{1}{t}\, dt\,+ \\[1.2em]
&\quad + \Biggl[ \frac{1}{t^{s}} Q(t)\Biggr]_1^\infty
  + s \int_1^\infty Q(t) t^{-s-1}\, dt
\end{aligned}
\end{equation}
by Stieltjes integration, which results in
\begin{equation}\label{eq:1}
Z(s) = \frac{1}{2\pi(s-1)^2} - \frac{\log(2\pi)}{2\pi(s-1)} + B_0 + s \int_{1}^{\infty} Q(t) t^{-s-1} \, dt
\end{equation}
where $B_0=-Q(1)$ as in (12). Since the integral in (17) converges for $\Re(s)>0$, it analytically continues $Z(s)$ in (17) to $\Re(s)>0$. It now remains to show that this integral admits the series expansion about $s=1$. This analysis is carried out by Bondarenko-Ivi$\acute{c}$-Saksman-Seip in [6, p.4-5] so we reproduce it here. We first consider an integral

\begin{equation}\label{eq:1}
\int_{1}^{\infty} t^{-s-1} Q(t)\, dt
\end{equation}
which converges for $\Re(s)>0$ assuming (RH). One then generates the exp-log expansion as
\begin{equation}\label{eq:1}
\int_{1}^{\infty} t^{-2} e^{-(s-1) \log t}Q(t) \, dt = \sum_{j=0}^{\infty} c_j \frac{(s-1)^j}{j!}
\end{equation}
where its expansion coefficients can be read off as
\begin{equation}\label{eq:1}
c_j = (-1)^j\int_{1}^{\infty} \frac{\log^j(t)}{t^{2}}Q(t)  \, dt
\end{equation}
but such a series is only limited to $|s-1|<1$ domain. However, on writing integral term in (17) in $s-1$ domain, then one needs to consider this form instead

\begin{equation}\label{eq:1}
\begin{aligned}
s\int_{1}^{\infty} t^{-s-1} Q(t)\, dt &= (s-1)\int_{1}^{\infty} t^{-s-1} Q(t)\, dt+\int_{1}^{\infty} t^{-s-1} Q(t)\, dt \\[1.2em]
&= \int_{1}^{\infty} [(s-1)t^{-s-1} Q(t)+t^{-s-1} Q(t)]\, dt \\[1.2em]
\end{aligned}
\end{equation}
where we can generate a new series expansion as
\begin{equation}\label{eq:1}
\begin{aligned}
s\int_{1}^{\infty} t^{-s-1} Q(t)\, dt&=\sum_{j=0}^{\infty}(j+1)c_j\frac{(s-1)^{j+1}}{(j+1)!}+\sum_{j=0}^{\infty}c_j\frac{(s-1)^{j}}{j!}\\[1.2em]
&=\sum_{m=1}^{\infty}mc_{m-1}\frac{(s-1)^{m}}{m!}+c_0+\sum_{m=1}^{\infty}c_m\frac{(s-1)^{m}}{m!}\\[1.2em]
&=c_0+\sum_{m=1}^{\infty}\left[mc_{m-1}+c_m\right]\frac{(s-1)^{m}}{m!}\\[1.2em]
&=\sum_{j=0}^{\infty} y_j \frac{(s-1)^j}{j!} \\[1.2em]
\end{aligned}
\end{equation}
by combining the expansions in (21). Here we've shifted index variable $m=j+1$ in first series in (22) and re-labeled $m=j$ in second series in (22). The new expansion coefficients are then

\begin{equation}\label{eq:1}
y_m=\begin{cases}
c_0 & \text{for } m = 0 \\
mc_{m-1}+c_m & \text{for } m \geq 1
\end{cases}
\end{equation}
where they are expressed by $c_n$ coefficients, and by (20) these coefficients can be expressed by the integral

\begin{equation}\label{eq:1}
y_m=(-1)^m\int_{1}^{\infty}\frac{\log^m(t)-m\log^{m-1}(t)}{t^2}Q(t)\, dt
\end{equation}
which matches the series in (15). And we note that for $m=0$ case we simplify (24) to get

\begin{equation}\label{eq:1}
y_0=\int_{1}^{\infty}\frac{Q(t)}{t^2}\, dt
\end{equation}
whereby adding the constant $B_0$ from (12) (and as it appears in (17)), we get the formula

\begin{equation}\label{eq:1}
C_0=-Q(1)+\int_{1}^{\infty}\frac{Q(t)}{t^2}\, dt
\end{equation}
And for the $m\geq 1$ case we have $C_m=y_m$. Also, we note the value of this constant
\begin{equation}\label{eq:1}
\begin{aligned}
-Q(1) &= L(1)\\[1.2em]
&=\frac{7}{8}+\frac{1}{2\pi}\log\left(\frac{1}{2\pi}\right)-\frac{1}{2\pi}\\[1.2em]
&= 0.42333783699382573900\ldots
\end{aligned}
\end{equation}
Hence this completes the proof.

\section{Numerical Computation}
\label{sec:comp}
In this section we explore numerical computation of this formula. To check the $n=0$ case, the formula (3) simplifies to

\begin{equation}\label{eq:1}
C_0 =\lim_{T\to\infty}  \Bigg\{\sum_{0<\gamma<T}\frac{1}{\gamma}-\frac{1}{4\pi}\log(T)\log\left(\frac{T}{4\pi^2}\right)\Bigg\}.
\end{equation}
This limit formula was first studied by Hassani [12][13] who first computed it accurately to $5$ decimal places. Hassani used a small number of non-trivial zeros and applied other approximation methods, such as averaging and curve fitting, to extract a few more digits. Later, Brent (BPT) computed this constant to $19$ decimal places [4, p. 64], using a database of $10^{10}$ non-trivial zeros together with (BPT) method, a new technique to improve convergence of sums over non-trivial zeros (Theorem 2 in [3]).  And in our previous article [16], with our best computation effort we computed

\begin{equation}\label{eq:20}
\begin{aligned}
C_0 = & 0.25163675131270596653346632934264537551475958738\\
      &3654550533059356530585960570182311791574050852516\\
      &937760994148142\ldots\\
\end{aligned}
\end{equation}
to $111$ decimal places using the Arias De Reyna algorithm [2]. But our result was offset by $\log^2(2\pi)/(4\pi)$ factor. And in Table $1$, we computed $C_n$ to $50$ digits for $0\leq n\leq 50$. We reproduce Table 1 here for reference.

So to verify equation (3) directly, we have a database of $2$ million non-trivial zeros from which we can easily compute such sums over non-trivial roots. We have $N(T)=2\times 10^6$ with

\begin{equation}\label{eq:20}
T= 1131944.47182486226849153321\ldots
\end{equation}
So to test the $n=0$ case, we directly compute the harmonic sum

\begin{equation}\label{eq:20}
\sum_{0<\gamma<T}\frac{1}{\gamma} = 11.63680321239784824535\ldots
\end{equation}
which results in calculation
\begin{equation}\label{eq:20}
C_0\approx 0.2516\textcolor{blue}{\underline{3}}729326778528275\ldots
\end{equation}
which is accurate to $5$ decimal places as it's shown by last correct digit underlined in blue. We did not employ any additional approximation techniques, as our goal is to just compute equation (3) directly. And similarly for $n=1$ case the formula yields

\begin{equation}\label{eq:1}
C_1 =\lim_{T\to\infty}  -\Bigg\{\sum_{0<\gamma<T}\frac{\log(\gamma)}{\gamma}-\frac{1}{12\pi}\log^2(T)\log\left(\frac{T^2}{8\pi^3}\right)\Bigg\}
\end{equation}
we compute the sum
\begin{equation}\label{eq:20}
\sum_{0<\gamma<T}\frac{\log(\gamma)}{\gamma} = 115.40475869669330537244\ldots
\end{equation}
which results in computing the constant
\begin{equation}\label{eq:20}
C_1\approx -0.130\textcolor{blue}{\underline{0}}05204046642585705\ldots
\end{equation}
which is accurate to $4$ decimal places as shown by last correct digit underlined in blue, as compared in the reference Table $1$. And similarly for $n=2$ we have

\begin{equation}\label{eq:1}
C_2 =\lim_{T\to\infty}  \Bigg\{\sum_{0<\gamma<T}\frac{\log^2(\gamma)}{\gamma}-\frac{1}{24\pi}\log^3(T)\log\left(\frac{T^3}{16\pi^4}\right)\Bigg\}
\end{equation}
and we compute the sum
\begin{equation}\label{eq:20}
\sum_{0<\gamma<T}\frac{\log^2(\gamma)}{\gamma} = 1238.24307097390760286254\ldots
\end{equation}
and the constant
\begin{equation}\label{eq:20}
C_2\approx 0.08\textcolor{blue}{\underline{2}}52679757366788286\ldots
\end{equation}
which matches to $3$ decimal places.

We note that as the accuracy starts dropping, and so more non-trivial roots are needed in computation to increase the accuracy.

\begin{table}[hbt!]
\caption{Reference table for $C_n$ coefficients computed by the (ADR) algorithm.}
\centering 
\begin{tabular}{| c | c |} 
\hline
$\boldsymbol{n}$  & $\boldsymbol{C_n}$ \textbf{(also in 50 digits)}\\ [0.5ex]
\hline
0 &  0.2516367513127059665334663293426453755147595873836 \\
\hline
1 & -0.1300444859118885707285274533988846777460553964263 \\
\hline
2 &  0.0824214912550528039526632284933172430791521350021 \\
\hline
3 & -0.0321581827282544905964296099391141952179545405019 \\
\hline
4 & -0.0531801364893419772868761573698112582469915802523  \\
\hline
5 &  0.2110321083617385257637243839874627961215847994456  \\
\hline
6 & -0.4933371057135871285817870279321636575675112589435  \\
\hline
7 &  0.9731261196976619662852108486791876458635644729040  \\
\hline
8 & -1.8021253179931622367536330625155209079039086674443  \\
\hline
9 &  3.7133510644596133576858937986178468541115390150895  \\
\hline
10& -11.583138616714443418004214394156033470878899508634  \\
\hline
\vdots & \vdots \\
\hline
20& -7.6931751083769270011123002218244304577221846239268$\times 10^{9}$ \\
\hline
30& -8.1910409909869137068367900925700302382658971132757$\times 10^{20}$ \\
\hline
40& -2.4605043425772457379890548734866774381567481629777$\times 10^{33}$ \\
\hline
50 &-8.9568228254793711194813512752380738598982095960590$\times 10^{46}$ \\
\hline
\end{tabular}
\label{table:nonlin} 
\end{table}

\section{Application of (BPT) for improving error bound in (3)}
In this Section, we will apply the work of Brent (BPT) Theorem to improve the error term in (14), which we rewrite as

\begin{equation}\label{eq:1}
\sum_{0<\gamma<T}\frac{\log^m(\gamma)}{\gamma}-A(T)=(-1)^mC_m+E(T)
\end{equation}
where

\begin{equation}\label{eq:1}
E(T) = O\left(\frac{\log^{m+1}(T)}{T}\right)
\end{equation}
as obtained earlier. This error is consistent with computations we did in Section $2$. Now applying (BPT) equations (1-3) [3, p.2] and taking $\phi(t)=\log^m(t)/t$, yields

\begin{equation}\label{eq:1}
E(T) = -\frac{\log^m(T)}{T}\Bigg[N(T)-L(T)\Bigg]+E_2(T)
\end{equation}
where the improved error term is

\begin{equation}\label{eq:1}
E_2(T) <  2(A_0+A_1\log T)\bigg\vert\frac{m\log^{m-1}(T)}{T^2}-\frac{\log^m (T)}{T^2}\bigg\vert+(A_1+A_2)\frac{\log^m (T)}{T^2}
\end{equation}
where the constants are $A_0= 2.067$, $A_1=0.059$ and $A_2= 0.007$. Hence

\begin{equation}\label{eq:1}
E_2(T) = O\left(\frac{\log^{m+1}(T)}{T^2}\right)
\end{equation}
a nice improvement by $O(1/T)$ factor from (40).

We summarize these results in the next Theorem:

\begin{theorem}
\label{thm:fta}
\begin{equation}\label{eq:1}
\begin{aligned}
C_n =\lim_{T\to\infty} &(-1)^n \Bigg\{\sum_{0<\gamma<T}\frac{\log^n(\gamma)}{\gamma}-\frac{1}{2\pi(n+1)(n+2)}\log^{n+1}(T)\log\left(\frac{T^{n+1}}{(2\pi)^{n+2}}\right)+\\
&-\frac{\log^n(T)}{T}\Bigg[N(T)-\frac{7}{8}-\frac{T}{2\pi}\log\left(\frac{T}{2\pi}\right)+\frac{T}{2\pi}\Bigg]+E_2(T)\Bigg\}
\end{aligned}
\end{equation}
for $n\geq 0$.
\end{theorem}

If we apply this Theorem to the parameters from the previous Section $2$ with $N(T)=2\times 10^{6}$ and $T$ as in (27), we recompute

\begin{equation}\label{eq:20}
C_0\approx 0.25163675131\textcolor{blue}{\underline{2}}53386814\ldots
\end{equation}
which is now accurate to $12$ decimal places as it's shown by last correct digit underlined in blue, and the computed error $E_2(T)< 4.5614\times 10^{-12}$. And similarly, recomputing
\begin{equation}\label{eq:20}
C_1\approx 0.130044485\textcolor{blue}{\underline{9}}0966041896\ldots
\end{equation}
which matches to $10$ decimal places, and the computed error $E_2(T)< 5.9073\times 10^{-11}$. And the next one is

\begin{equation}\label{eq:20}
C_2\approx 0.082421491\textcolor{blue}{\underline{2}}2637581623 \ldots
\end{equation}
which matches to $10$ decimal places, and the computed error $E_2(T)< 7.6058\times 10^{-10}$.

This demonstrates that by adding a simple extra term from (BPT), we can get a big improvement in convergence. We also numerically verified this formula for higher $n$, and it is valid as expected.

\texttt{Email: art.kawalec@gmail.com}

\end{document}